\documentclass[12pt]{article}
\usepackage{hyperref}
\usepackage{amsfonts}
\usepackage{enumerate}
\usepackage{amsthm}
\usepackage{amsmath}
\usepackage{graphicx}

\begin{document}

\begin{center}
{\large\bf On the Unity Row Summation and Real Valued Nature of the $F_{LG}$ Matrix}

\vskip.20in

Ioannis K. Dassios$^{1, 2
}$, Paul Cuffe$^{2}$, Andrew Keane$^{2}$ \\[2mm]
{\footnotesize
$^{1}$MACSI, Department of Mathematics \& Statistics, University of Limerick, Ireland\\[5pt]
$^{2}$Electricity Research Centre, University College Dublin, Ireland\\[5pt]
}
\end{center}

{\footnotesize
\noindent
\textbf{Abstract.} Electrical power system calculations rely heavily on the $Y_{bus}$ matrix, which is the Laplacian matrix of the network under study, weighted by the complex-valued admittance of each branch. It is often useful to partition the $Y_{bus}$ into four submatrices, to separately quantify the connectivity between and among the load and generation nodes in the network. Simple manipulation of these submatrices gives the $F_{LG}$ matrix, which offers useful insights on how voltage deviations propagate through a power system and on how energy losses may be minimized. Various authors have observed that in practice the elements of $F_{LG}$ are real-valued and its rows sum close to one: the present paper explains and proves these properties. \\[3pt]
\\[3pt]

\vskip.2in

\section{ Introduction}

Currents $(I)$ and voltages $(V)$ in an electrical power system are related by the admittance matrix, $Y_{bus}$, which is generally constructed to have the properties of a weighted Laplacian matrix (disregarding shunt element modelling), see \cite{bib1}, \cite{bib2}. It can usefully be partitioned as follows by reordering to group generator $(G)$ and load $(L)$ nodes separately:
\begin{equation}\label{eq1}
\left[\begin{array}{c}I_G\\I_L\end{array}\right]=Y_{bus}\left[\begin{array}{c}V_G\\V_L\end{array}\right].
\end{equation}
Where $I_G, V_G\in\mathbb{C}^m$, $I_L, V_L\in\mathbb{C}^n$ and 
\[
Y_{bus}=\left[\begin{array}{cc}Y_{GG}&Y_{GL}\\Y_{LG}&Y_{LL}\end{array}\right].
\]
Where $Y_{GG}\in\mathbb{C}^{m\times m}$, $Y_{LL}\in\mathbb{C}^{n\times n}$ and $Y_{LG}=Y_{GL}^T\in\mathbb{C}^{n\times m}$. With $()^T$ we denote the (non-conjugate) transposed tensor. The $Y_{bus}$ matrix is square by definition. It has $N$ rows and $N$ columns, where $N$ is the number of nodes in the power system. The system is reordered and partitioned to group $G$ buses and $L$ buses separately. The $n$ $L$ buses plus the $m$ $G$ buses equals the total numbers of buses $N$. Although there may be buses in a power system that don't connect any $L$ or $G$ but are just passive interconnection points, in this analysis they would be grouped with the load buses. Typically $Y_{LL}$ will be larger. For instance, on the $IEEE$ 118 bus system we have 19 $G$ nodes and 99 $L$ nodes, i.e. $Y_{GG}\in\mathbb{C}^{19\times 19}$ and $Y_{LL}\in\mathbb{C}^{99\times 99}$. The other sub-matrices are not generally square, in this case $Y_{LG}\in\mathbb{C}^{99\times 19}$ and $Y_{GL}\in\mathbb{C}^{19\times 99}$. Note that $Y_{LG}=Y^T_{GL}$. The justification for this assertion is the guaranteed symmetry of the $Y_{bus}$ matrix (it is constructed to have this property, which is only broken under the rare circumstance that phase-shifting transformers are included in the model).
The $Y_{bus}$ matrix is complex-valued. The physical properties of electrical conductors imply that the imaginary part of each element is larger than the real part in nearly all cases. In the below figure is plotted the real and imaginary component of every element of the $Y_{bus}$ matrix for the 118 bus system. The imaginary component is on the vertical axis:
\begin{figure}[h]
\begin{center}
\includegraphics[width=0.50\textwidth]{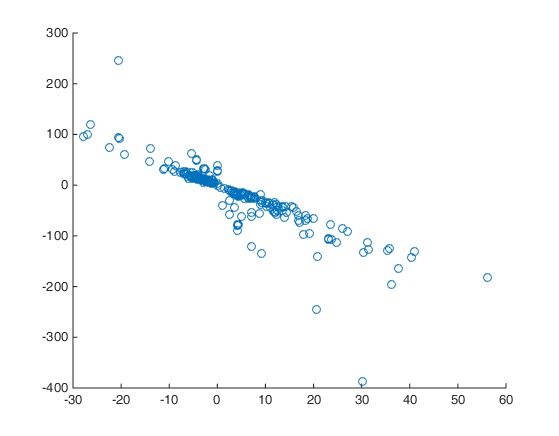}
\caption{Empirical Real-Imaginary ratio}
\label{figure1}
\end{center}
\end{figure}
The real:imaginary ratio is approximately equal for all elements, as can be seen in Figure 1.

The partitioning of the $Y_{bus}$ matrix into load and generation blocks was introduced by Kessel and Glavitsch in \cite{bib3}, and has subsequently been applied to a wide range of power engineering problems  \cite{bib4}, \cite{bib5}, \cite{bib6}, \cite{bib7}, \cite{bib8}, \cite{bib9}, \cite{bib10}, \cite{bib11}, \cite{bib12}. 
Recent work by Sikiru et. al. \cite{bib13}, \cite{bib14}, \cite{bib15}  has used the partitioned $Y_{bus}$ approach to develop a more fundamental understanding of a power systems inherent connective structure, using, for instance, Schur complements and eigen analysis. Relatedly, recent work by Abdelkader et al.  \cite{bib16}, \cite{bib17}, \cite{bib18} offers a deeper conceptual understanding of these $Y_{bus}$ partitions, demonstrating how they allow the separation of currents in the network into load and generator induced components. Notably \cite{bib18} demonstrates how a strictly equal real:imaginary ratio for every element in the $Y_{bus}$ matrix brings one component of the network physical power losses to zero.
 
From \eqref{eq1}:
\begin{equation}\label{eq2}
I_G=Y_{GG}V_G+Y_{GL}V_L
\end{equation}
and 
\begin{equation}\label{eq3}
I_L=Y_{LG}V_G+Y_{LL}V_L.
\end{equation}
Rearranging \eqref{eq3}:
\begin{equation}\label{eq4}
V_L=Z_{LL}I_L-Z_{LL}Y_{LG}V_G.
\end{equation}
Where
\[
Z_{LL}=\left\{\begin{array}{cc} Y_{LL}^{-1},&if\quad detY_{LL}\neq 0\\Y_{LL}^{\dagger},&if\quad det(Y_{LL})=0\end{array}\right\}.
\]
The matrix $Y_{LL}^\dagger$ is the Moore-Penrose Pseudoinverse of $Y_{LL}$, calculated by the singular value decomposition of $Y_{LL}$, see \cite{bib19}, \cite{bib20}. Substituting for $V_L$ in \eqref{eq2}:
\begin{equation}\label{eq5}
I_G=(Y_{GG}-Y_{GL}Z_{LL}Y_{LG})V_G+Y_{GL}Z_{LL}I_L.
\end{equation}
Equations \eqref{eq4} and \eqref{eq5} are typically represented in matrix form, which permits useful engineering applications:
\[
\left[\begin{array}{c}V_L\\I_G\end{array}\right]=\left[\begin{array}{cc}Z_{LL}&F_{LG}\\K_{GL}&Y_{GGM}\end{array}\right]\left[\begin{array}{c}I_L\\V_G\end{array}\right].
\]
Where:
\[
Y_{GGM}=Y_{GG}-Y_{GL}Z_{LL}Y_{LG}
\]
and
\begin{equation}\label{eq8}
F_{LG}=-Z_{LL}Y_{LG}=-K_{GL}^T
\end{equation}
Various works (e.g \cite{bib14}, \cite{bib4}, \cite{bib21}) have noted, and, indeed, relied upon, the observation that the elements of $F_{LG}$ are in practice real-valued and that its rows sum to unity, or close to. 

\section{Main results}
As written in the previous section, the ratio $(Im:Re)$ of each entry of the $Y_{bus}$ matrix will in practice tend to be fairly homogeneous though this isn't always guaranteed. The $Y_{bus}$ matrix will have generally diagonal elements that are positive real and negative imaginary. Off-diagonals will be negative real and positive imaginary. Based on these observations we can provide the following Proposition.
\\\\
\textbf{Proposition 2.1.} Assume the $Y_{bus}$ matrix as defined in \eqref{eq1}. If the entries in each row have the real:imaginary ratio equal, then all entries of the $F_{LG}$ matrix, defined in \eqref{eq8}, are real numbers. In addition if the $Y_{bus}$ matrix has diagonal elements that are non-negative real, non-positive imaginary and off diagonals are non-positive real, non-negative imaginary, then all entries of $F_{LG}$, are non-negative real numbers.
\\\\
\textbf{Proof.} Let
\[
Y_{LG}=  \left[\begin{array}{cccc}a_{11}  & a_{12} & \dots  & a_{1m}  \\
   a_{21} & a_{22} &   \dots&a_{2m}  \\
    \vdots  &  \vdots  &  \ddots  &   \vdots   \\
   a_{n1} & a_{n2} &  \dots& a_{nm}
   \end{array}\right],\quad Y_{LL}=\left[\begin{array}{cccc} b_{11}  & b_{12} & \dots  & b_{1n}  \\b_{21} & b_{22} &   \dots&b_{2n}  \\
    \vdots  &  \vdots  &  \ddots  &   \vdots   \\
   b_{n1} & b_{n2} &  \dots& b_{nn}
   \end{array}\right]
\]
and $a_{kj}=Re(a_{kj})+iIm(a_{kj})$, $b_{kj}=Re(b_{kj})+iIm(b_{kj})$. Then for each row $k=1,2,...,n$ of $Y_{LG}$, $Y_{LL}$ we have
\[
\frac{Im(a_{kj})}{Re(a_{kj})}=\frac{Im(b_{kj})}{Re(b_{kj})}=u_k\in\mathbb{R},\quad k=1,2,...,n.
\]
From \eqref{eq8} we have
\[
Y_{LL}F_{LG}=-Y_{LG},
\]
or, equivalently,
\[
\left[\begin{array}{cccc} b_{11}  & b_{12} & \dots  & b_{1n}  \\b_{21} & b_{22} &   \dots&b_{2n}  \\
    \vdots  &  \vdots  &  \ddots  &   \vdots   \\
   b_{n1} & b_{n2} &  \dots& b_{nn}
   \end{array}\right]F_{LG}=
   -\left[\begin{array}{cccc}a_{11}  & a_{12} & \dots  & a_{1m}  \\
   a_{21} & a_{22} &   \dots&a_{2m}  \\
    \vdots  &  \vdots  &  \ddots  &   \vdots   \\
   a_{n1} & a_{n2} &  \dots& a_{nm}
   \end{array}\right],
  \] 
  or, equivalently,
\[
\left[\begin{array}{cccc} Re(b_{11})[1+iu_1]  & Re(b_{12})[1+iu_1] & \dots  & Re(b_{1n})[1+iu_1]  \\Re(b_{21})[1+iu_2] & Re(b_{22})[1+iu_2] &   \dots&Re(b_{2n}) [1+iu_2] \\
    \vdots  &  \vdots  &  \ddots  &   \vdots   \\
   Re(b_{n1})[1+iu_n] & Re(b_{n2})[1+iu_n] &  \dots& Re(b_{nn})[1+iu_n]
   \end{array}\right]F_{LG}=
   \]
   \[-\left[\begin{array}{cccc}Re(a_{11})[1+iu_1]  & Re(a_{12})[1+iu_1] & \dots  & Re(a_{1m}) [1+iu_1] \\
   Re(a_{21}) [1+iu_2]& Re(a_{22}) [1+iu_2]&   \dots&Re(a_{2m})[1+iu_2]  \\
    \vdots  &  \vdots  &  \ddots  &   \vdots   \\
   Re(a_{n1})[1+iu_n] & Re(a_{n2}) [1+iu_n]&  \dots& Re(a_{nm})[1+iu_n]
   \end{array}\right],
  \]   
or, equivalently,
\[
diag\left\{\begin{array}{cccc}1+iu_1,& 1+iu_2,& \dots,& 1+iu_n\end{array}\right\} \left[\begin{array}{cccc} Re(b_{11})  & Re(b_{12}) & \dots  & Re(b_{1n})  \\Re(b_{21}) & Re(b_{22})&   \dots&Re(b_{2n}) \\
    \vdots  &  \vdots  &  \ddots  &   \vdots   \\
   Re(b_{n1}) & Re(b_{n2}) &  \dots& Re(b_{nn})
   \end{array}\right]F_{LG}=
   \]
   \[
   -diag\left\{\begin{array}{cccc}1+iu_1,& 1+iu_2,& \dots,& 1+iu_n\end{array}\right\}\left[\begin{array}{cccc}Re(a_{11}) & Re(a_{12}) & \dots  & Re(a_{1m}) \\
   Re(a_{21}) & Re(a_{22})&   \dots&Re(a_{2m})  \\
    \vdots  &  \vdots  &  \ddots  &   \vdots   \\
   Re(a_{n1}) & Re(a_{n2}) &  \dots& Re(a_{nm})
   \end{array}\right].
  \]
 or, equivalently, since $det(diag\left\{\begin{array}{cccc}1+iu_1,& 1+iu_2,& \dots,& 1+iu_n\end{array}\right\})\neq 0$
   \[
Re(Y_{LL})F_{LG}=-Re(Y_{LG}).
 \]  
 This is a liner system where the known matrices have only real entries. Hence, any solution for $F_{LG}$ lies inside of $\mathbb{R}^{n\times m}$, i.e. all entries of $F_{LG}$ are real numbers. In addition since the $Y_{bus}$ matrix is assumed to have diagonal elements that are non-negative real \& non-positive imaginary and off diagonals are non-positive real \& non-negative imaginary, the matrix $Re(Y_{LG})$ has all its entries non-positive, i.e. the matrix $Re(Y_{LG})$ has all its entries non-negative and the matrix $Re(Y_{LL})$ has its diagonal elements non-negative and its off diagonals non-positive, i.e. the inverse (or pseudo inverse) of $Re(Y_{LL})$ has all its entries non-negative, see \cite{bib22}. Hence the matrix $F_{LG}$ will have all its elements non-negative since it is given from the product of two matrices, the inverse (or pseudo inverse) of $Re(Y_{LL})$ and the matrix $Re(Y_{LG})$ . The proof is completed. 
 \\\\
If we use a short line model of the system and so neglect shunt elements, then each row of the $Y_{bus}$ matrix sums to zero. If we include shunts, then each row sums close to zero. We can state the following Theorem.\\\\
\textbf{Theorem 2.1.} Assume the $Y_{bus}$ matrix as defined in \eqref{eq1} and the matrix $F_{LG}$, as defined in \eqref{eq8}. Then the row sum of $F_{LG}$ is one if the rows of $Y_{bus}$ sum to zero and $det(Y_{LL})\neq 0$ and is close to one if the rows of $Y_{bus}$ sum close to zero or sum to zero and $det(Y_{LL})=0$.
\\\\
\textbf{Proof.} Let for system \eqref{eq1}
\begin{equation}\label{eqa}
Y_{LG}=[a_{ij}]^{j=1,2,...,m}_{i=1,2,...,n}
\end{equation}
and
\begin{equation}\label{eqb}
Y_{LL}=[b_{ij}]^{j=1,2,...,n}_{i=1,2,...,n}.
\end{equation}
Where
\[
[a_{ij}]^{j=1,2,...,m}_{i=1,2,...,n}=  \left[\begin{array}{cccc}a_{11}  & a_{12} & \dots  & a_{1m}  \\
   a_{21} & a_{22} &   \dots&a_{2m}  \\
    \vdots  &  \vdots  &  \ddots  &   \vdots   \\
   a_{n1} & a_{n2} &  \dots& a_{nm}
   \end{array}\right],\quad [b_{ij}]^{j=1,2,...,n}_{i=1,2,...,n}=\left[\begin{array}{cccc} b_{11}  & b_{12} & \dots  & b_{1n}  \\b_{21} & b_{22} &   \dots&b_{2n}  \\
    \vdots  &  \vdots  &  \ddots  &   \vdots   \\
   b_{n1} & b_{n2} &  \dots& b_{nn}
   \end{array}\right].
\]
We assume that each row of $Y_{bus}$ sums to zero, i.e. $\forall$ row $i=1, 2, ..., n$ of the matrices \eqref{eqa}, \eqref{eqb} we have $\sum^m_{j=1}a_{ij}+\sum^n_{j=1}b_{ij}=0$, or equivalently
\begin{equation}\label{eqc}
\sum^m_{j=1}a_{ij}=-\sum^n_{j=1}b_{ij},\quad \forall i=1,2,...,n.
\end{equation}
From \eqref{eq8} we have
\begin{equation}\label{eqd}
Y_{LL}F_{LG}=-Y_{LG},
\end{equation}
Let
\begin{equation}\label{eqe}
F_{LG}=[c_{ij}]^{j=1,2,...,m}_{i=1,2,...,n}.
\end{equation}
Then $\forall i=1,2,...,n$ by substituting \eqref{eqa}, \eqref{eqb}, \eqref{eqe} into \eqref{eqd} and by calculating every component in row $i$, we get
\[
\begin{array}{c}
     -a_{i1}=\sum_{k=1}^n b_{ik}c_{k1} \\
     -a_{i2}=\sum_{k=1}^n b_{ik}c_{k2}\\
     \vdots\\
     -a_{im}=\sum_{k=1}^n b_{ik}c_{km}.
    \end{array}
\]
By taking the sum of the above equalities, $\forall i=1,2,...,n$ we arrive at
\[
-[a_{i1}+a_{i2}+..+a_{im}]=\sum_{k=1}^n b_{ik}c_{k1}+\sum_{k=1}^n b_{ik}c_{k2}+\dots+\sum_{k=1}^n b_{ik}c_{km},
\]
or, equivalently,
\[
-\sum_{j=1}^ma_{ij}=\sum_{j=1}^m\sum_{k=1}^n b_{ik}c_{kj}.
\]
By using \eqref{eqc} on the above expression we get
\[
\sum_{j=1}^nb_{ij}=\sum_{j=1}^m\sum_{k=1}^n b_{ik}c_{kj},
\]
or, equivalently,
\[
\sum_{k=1}^nb_{ik}=\sum_{j=1}^m\sum_{k=1}^n b_{ik}c_{kj},
\]
or, equivalently,
\[
b_{i1}+b_{i2}+\dots+b_{in}=\sum_{j=1}^m[b_{i1}c_{1j}+b_{i2}c_{2j}+\dots+b_{in}c_{nj}],
\]
or equivalently,
\[
\begin{array}{cc}b_{i1}+b_{i2}+\dots+b_{in}=&b_{i1}c_{11}+b_{i2}c_{21}+\dots+b_{in}c_{n1}+\\
&b_{i1}c_{12}+b_{i2}c_{22}+\dots+b_{in}c_{n2}+\\
&\vdots\\
&b_{i1}c_{1m}+b_{i2}c_{2m}+\dots+b_{in}c_{nm},\end{array}
\]
or, equivalently,
\[
\begin{array}{cc}b_{i1}+b_{i2}+\dots+b_{in}=&b_{i1}(c_{11}+c_{12}+\dots+c_{1m)}+\\
&b_{i2}(c_{21}+c_{22}+\dots+c_{2m})+\\
&\vdots\\
&b_{in}(c_{n1}+c_{n2}+\dots+c_{nm}),\end{array}
\]
or, equivalently,
\[
b_{i1}+b_{i2}+\dots+b_{in}=b_{i1}(\sum_{j=1}^mc_{1j})+b_{i2}(\sum_{j=1}^mc_{2j})+\dots+b_{in}(\sum_{j=1}^mc_{nj}),
\]
or, equivalently,
\begin{equation}\label{eqf}
b_{i1}[(\sum_{j=1}^mc_{1j})-1]+b_{i2}[(\sum_{j=1}^mc_{2j})-1]+\dots+b_{in}[(\sum_{j=1}^mc_{nj})-1]=0.
\end{equation}
Let $\forall i=1,2,...,n$ 
\begin{equation}\label{eqg}
q_i=(\sum_{j=1}^mc_{ij})-1.
\end{equation}
Then by replacing \eqref{eqg} into \eqref{eqf}, $\forall i=1,2,...,n$ we have
\[
b_{i1}q_1+b_{i2}q_2+\dots+b_{in}q_n=0,
\]
or, equivalently,
\[
\left[\begin{array}{cccc}b_{i1}&b_{i2}&\dots&b_{in}\end{array}\right]\left[\begin{array}{c}q_1\\q_2\\\vdots\\q_n\end{array}\right]=0,\quad \forall i=1,2,...,n.
\]
Since the above expression holds $\forall i=1,2,...,n$ we have equivalently
\[
\left[\begin{array}{cccc}b_{11}  & b_{12} & \dots  & b_{1n}  \\
   b_{21} & b_{22} &   \dots&b_{2n}  \\
    \vdots  &  \vdots  &  \ddots  &    \vdots   \\
   b_{n1} & b_{n2} &  \dots& b_{nn}
   \end{array}\right]\left[\begin{array}{c}q_1\\q_2\\\vdots\\q_n\end{array}\right]=\left[\begin{array}{c}0\\0\\\vdots\\0\end{array}\right].
\]
By setting $q=\left[\begin{array}{c}q_1\\\vdots\\q_n\end{array}\right]$, $0_{n,1}=\left[\begin{array}{c}0\\0\\\vdots\\0\end{array}\right]$ and by using \eqref{eqb} we have
\[
Y_{LL}\cdot q=0_{n,1}.
\]
From the above expression if $det(Y_{LL})\neq0$, then $q=0_{n,1}$, or, equivalently, $q_i=0,$ $\forall i=1,2,...,n$, i.e. by using \eqref{eqg}
\begin{equation}\label{eqh}
\sum_{j=1}^mc_{ij}=1.
\end{equation}
If $det(Y_{LL})=0$, then $q\cong 0_{n,1}$ (by using the pseudo inverse of $Y_{LL}$ via the SVD method), or, equivalently, $q_i\cong0,$ $\forall i=1,2,...,n$, i.e. by using \eqref{eqg}
\begin{equation}\label{eqi}
\sum_{j=1}^mc_{ij}\cong1.
\end{equation}
Thus, from \eqref{eqh} and \eqref{eqi} every row of $F_{LG}$ will sum to 1 if $det(Y_{LL})\neq 0$ and will sum closely to 1 if $det(Y_{LL})=0$. If the rows of $Y_{bus}$ sum close to zero, then with similar steps we arrive at \eqref{eqi}. The proof is completed.

\section*{Conclusions}
The partitioning of the $Y_{bus}$ matrix has opened numerous fruitful avenues in power system analysis, and it is hoped that a clearer understanding of the matrix properties underpinning these partitions will support this strand of research. This work has proved two numerically observed properties of the $F_{LG}$ matrix: future work may consider the matrix characteristics of the other sub-matrices derived by the $Y_{bus}$ partitioning approach.

\subsection*{Acknowledgments}
I. Dassios is supported by Science Foundation Ireland (award 09/SRC/E1780). P. Cuffe is funded through the European Community's Seventh Framework Programme (FP7/2007-2013) under grant agreement nA$^o$608732A. This work was conducted in the Electricity Research Centre, University College Dublin, Ireland, which is supported by the Commission for Energy Regulation, Bord Gais Energy, Bord na Mona Energy, Cylon Controls, EirGrid, Electric Ireland, Energia, EPRI, ESB International, ESB Networks, Gaelectric, Intel, SSE Renewables, and UTRC.

\end{document}